\documentclass[12pt]{article}

\usepackage{amsthm}

\usepackage[left=2cm,right=2cm,top=2cm,bottom=2cm,bindingoffset=0cm]{geometry}

\usepackage[T2A]{fontenc}
\usepackage[utf8]{inputenc}
\usepackage[english]{babel}
\usepackage{amsmath}
\usepackage{amsfonts}
\usepackage{mathtools}
\usepackage{graphicx}
\usepackage{systeme}
\usepackage{tikz}
\usepackage{blindtext}
\usepackage{amssymb}
\usepackage{hyperref}
\usepackage{bbm}
\usepackage{cite}
\usepackage[normalem]{ulem}
\graphicspath{ {./} }

\newtheorem{theorem}{Theorem}[section]

\newtheorem{lemma}{Lemma}[section]
\newtheorem*{lemma*}{Lemma}
\newtheorem*{den*}{Denotion}

\newtheorem*{remark*}{Remark}

\theoremstyle{remark}
\newtheorem*{remark}{Remark}

\theoremstyle{definition}
\newtheorem{definition}{Definition}[section]

\title{Spectral gap of $G(n, \alpha n, \alpha^2 n)$ graphs and the giant component theorem\footnote{The research is supported by the Ministry of Science and Higher Education of the Russian Federation, project No. FSMG-2024-0011 and the grant of Theoretical Physics and Mathematics Advancement Foundation “BASIS”.}}
\author{M. Koshelev\footnote{Moscow Institute of Physics and Technology, Moscow State University, mkoshelev99@gmail.com}, A. Raigorodskii\footnote{Moscow Institute of Physics and Technology, Moscow State University, Institute of Mathematics and Computer Science, Buryat State University; Caucasus Mathematical Centre, Adyghe State University, mraigor@yandex.ru}}
\date{}

\begin{document}

\maketitle

\begin{abstract}
    The spectrum of a graph $G$ is the set of the eigenvalues of its adjacency matrix. It turns out that one can say a lot about a graph with the only knowledge being the spectrum of this graph. In this paper we obtain new results about the spectrum of $G(n, \alpha n, \alpha^2 n)$ graphs. We then apply these results to get a giant component theorem for them.
\end{abstract}

\section{Introduction}

In this paper we mainly work with Johnson graphs, also known as $G(n, r, s)$ graphs. The vertices of such graphs are $n$-dimensional vectors $v \in \{0,1\}^n$ such that $\|v\|^2 = r$. An edge is drawn between two vertices if and only if the inner product of the corresponding vectors is equal to $s$.

Johnson graphs have found their applications in various topics of combinatorial geometry, Ramsey theory, coding theory and other branches of modern combinatorics. There exists a vast literature on such graphs, see, for example,~\cite{Codes, KK, FrW, LS, RaiMeta, RaiMetaNew}.

The property we consider in this paper is the spectrum of such graphs. It has received a lot of attention in the past 50 years. It appears that one can say a lot about a graph just by looking at its eigenvalues. Before we continue, let us give some definitions.

\begin{definition}
    Let $G$ be a $d$-regular graph on $n$ vertices. Let $A$ be the adjacency matrix of $G$. We then call the sequence $\lambda_1 \geq \lambda_2 \dotsc \geq \lambda_n$ of the eigenvalues of $A$ the spectrum of $G$. We denote the degree of a vertex in a $d$-regular graph $G$ as $d(G)$.
\end{definition}

\begin{definition}
Let $G$ be a $d$-regular graph on $n$ vertices. Then 
    $$
        \lambda(G) := \max\{\lambda_2, -\lambda_n\} = \max\{|\lambda_2|, \dotsc, |\lambda_n|\}.
    $$
\end{definition}

The following theorem about the spectrum of $G(n, r, s)$ graphs is known (see \cite{Spectrum}, Theorem 4.6).

\begin{theorem}[\cite{Spectrum}]
\label{Eigenvalues}
The set of eigenvalues of $G(n, r, s), r \leq n/2$ coincides with the set $\{E_{r - s}(i), 0 \leq i \leq r\}$, where 
\begin{multline*}
    E_k(i) = \sum\limits_{j=0}^{k}{(-1)^{k-j}\binom{r-j}{k-j}\binom{r-i}{j}\binom{n-r+j-i}{j}} = \sum\limits_{j = 0}^{k}{(-1)^{j}\binom{i}{j}\binom{r - i}{k - j}\binom{n - r - i}{k - j}}.
\end{multline*}
Moreover, the multiplicity of $ E_{r - s}(i)$ is equal to $\binom{n}{i} - \binom{n}{i - 1}$. Note that in these formulas some of the binomial coefficients $\binom{n}{k}$ might have either $k < 0$ or $k > n$. In this case we assume them to be equal to zero.
\end{theorem}


Unfortunately, we cannot deduce the value of $\lambda(G(n,r,s))$ immediately from this theorem, as it is not obvious how to compare the values of $E_k(i)$ for different $i$. This means some additional work needs to be done to obtain the value of $\lambda(G(n,r,s))$. There are a few results of this kind.

The first one is due to Lov\'asz and it gives the value of $\lambda(G(n,r,0))$, also known as Kneser graphs.

\begin{theorem}[\cite{LovLambda}]
\label{LovT}
$\lambda(G(n, r, 0)) = \frac{r}{n - r}d(G(n, r, 0))$.
\end{theorem}

The second one was introduced in \cite{Brouwer} and works for a wider range of parameters.

\begin{theorem}[\cite{Brouwer}]
\label{BrT}
Let $(r - s)(n - 1) \geq r(n - r)$. Then $$\lambda(G(n, r, s)) = |E_{r - s}(1)| = \frac{|sn - r^2|}{r(n - r)}d(G(n, r, s)).$$
\end{theorem}

Some results of this kind were also obtained in \cite{KosEig}. 

The first one of these results covers one special case of parameters, namely, $n = 2r = 4s$.

\begin{theorem}
\label{Th2}
There exists an integer $N$ such that for all $s > N$
$$
\lambda(G(4s, 2s, s)) = |E_{s}(2)| = \frac{4s - 2}{s^2}\binom{2s - 2}{s - 1}^2 = \frac{1}{4s - 2}E_s(0) = O\left(\frac{d(G(4s, 2s, s))}{s}\right).
$$
\end{theorem}

The second one covers a wide range of parameters at the cost of providing only some bounds on $\lambda$ instead of the exact values.

\begin{theorem}
\label{ThMerge}
\begin{enumerate}
    \item Let $r > s \geq 1$. Then there exists an integer $N$ such that $\lambda(G(n, r, s)) = E_{r - s}(1) \sim \frac{s}{r}d(G(n, r, s))$ for all $n > N$.
    \item Let $n \gg r(n) \gg s(n)$, $r^2(n) \gg n$. Then $\lambda(G(n, r, s)) = o(d(G(n, r, s))).$
    \item Let $r(n) = O(\sqrt{n}), r \gg s \geq 1$. Then $\lambda(G(n, r, s)) = O\left(\frac{s}{r}d(G(n, r, s))\right).$
    Moreover, if at least one of  
    \begin{enumerate}
    \item $s = \omega(1)$;
    \item $r = o(\sqrt{n})$,
    \end{enumerate}
    holds, then $\lambda(G(n, r, s)) = (1 + o(1))\frac{s}{r}d(G(n, r, s)).$
\end{enumerate}
\end{theorem}

\begin{remark}
In this text we use the denotation $f(n) \gg (\ll) g(n)$ for functions $f$ and $g$ such that $\frac{f(n)}{g(n)} \to \infty (\to 0)$ as $n \to \infty$.
\end{remark}

In this paper we will try to generalize Theorem \ref{Th2} for $G(n, \alpha n, \alpha^2 n)$ graphs and also provide one possible application of the obtained results.

The rest of the paper is organized in the following way. In Section 2 we state the main result. Section 3 contains the proof of this result. Finally, in Section 4 we apply the main theorem to obtain the result about the emergence of giant component in $G(n,r,s)$ graphs. 

\section{Formulation of the results}

There are (at least) two natural ways to generalize Theorem \ref{Th2}. For example, one may try to replace the values 4 and 2 in $n = 2r = 4s$ constraint with something else. The other idea is to prove a similar result not only for $G(4s,2s,s)$, but also for $G(4s + f_1(s),2s + f_2(s),s + f_3(s))$, where $f_i(s)$ are small in some sense.

In this paper we will do the generalization in both senses simultaneously at the cost of moving from calculating the exact value of $\lambda$ to only bounding it from above. The main theorem of this paper is the following.

\begin{theorem}
    \label{Main}
    Let us fix $\alpha \in (0,1)$ and two functions $f_r(n)$ and $f_s(n)$, $f_s(n), f_r(n) = o\left(n\right)$. Then the following equality holds.

    \begin{multline*}
        \lambda\left(G\left(n,\alpha n + f_r(n), \alpha^2n + f_s(n)\right)\right) = \\ = O\left(\frac{\max\{1, |f_s(n) - 2\alpha f_r(n)|, f_r^2(n)/n\}}{n}d\left(G\left(n,\alpha n + f_r(n), \alpha^2n + f_s(n)\right)\right)\right).
    \end{multline*}
\end{theorem}

\section{Proof}

In the proof we will use the following lemma from \cite{KosEig}.

\begin{lemma}
\label{exp2}
Let $C_{n, r, s}(i, j) = \binom{i}{j}\binom{r - i}{r - s - j}\binom{n - r - i}{r - s - j}$. Then 
\begin{multline*}
    C_{n, r, s}(i, j) = \frac{(r-s-j+1)^2}{(r - i + 1)(n - r - i + 1)}C_{n, r, s}(i - 1, j - 1) +\\+ \frac{(s - i + j + 1)(n - 2r + s + j - i + 1)}{(r - i + 1)(n - r - i + 1)}C_{n, r, s}(i - 1, j).
\end{multline*}
\end{lemma}

\begin{proof}[Proof of Theorem \ref{Main}]
    Define $r = \alpha n + f_r(n), s = \alpha^2 n + f_s(n)$. We may assume that either $\alpha < 1/2$ or $\alpha = 1/2$ and $f_r(n) \leq 0$. Indeed, suppose that $\alpha > 1/2$. Let $\alpha^{\prime} = 1 - \alpha$, $f_r^{\prime}(n) = -f_r(n)$, $f_s^{\prime} = -2f_r(n)+f_s(n)$. It is easy to see that $(f_r^{\prime})^2(n) = f_r^2(n)$. $f_s^{\prime}(n) - 2\alpha^{\prime}f_r^{\prime}(n) = f_s(n) - 2\alpha f_r(n)$ and, finally, $G(n,n-r,n-2r+s) = G\left(n,\alpha^{\prime} n + f_r^{\prime}(n), (\alpha^{\prime})^2n + f_s^{\prime}(n)\right)$. It remains to notice that the graphs $G(n,r,s)$ and $G(n,n-r,n-2r+s)$ are isomorphic. The same trick applies to the case where $\alpha = 1/2$ and $f_r(n) > 0$. From now on we will assume that $n \geq 2r$.
    
    Using Lemma \ref{exp2}, we obtain (for $i \leq r - s$)
    \begin{multline*}
    \label{eq1}
        E_{r-s}(i) = \sum_{j = 0}^{i}(-1)^jC_{n,r,s}(i,j) = \\ =  \sum_{j = 0}^{i}(-1)^j\left(\frac{(r-s-j+1)^2}{(r - i + 1)(n - r - i + 1)}C_{n, r, s}(i - 1, j - 1) +\right.\\\left.+ \frac{(s - i + j + 1)(n - 2r + s + j - i + 1)}{(r - i + 1)(n - r - i + 1)}C_{n, r, s}(i - 1, j)\right) = \\ = \sum_{j = 0}^{i-1}(-1)^j\frac{(s - i + j + 1)(n - 2r + s + j - i + 1) - (r-s-j)^2}{(r - i + 1)(n - r - i + 1)}C_{n,r,s}(i-1,j) = \\ = \sum_{j = 0}^{i-1}(-1)^j\frac{i^2 - 2 i j - i n + 2 i r - 2 i s - 2 i + j n + 2 j + n s + n - r^2 - 2 r + 2 s + 1}{(r - i + 1)(n - r - i + 1)}C_{n,r,s}(i-1,j)
    \end{multline*}

    The next thing to do is to obtain some bound on $C_{n,r,s}(i,j)$. Let 
    $$\Tilde{C}_{n,r,s}(i) = \max\limits_{j \in 0, \dotsc, \min(r-s, i)}C_{n,r,s}(i, j).$$

    We will now prove that for all $\varepsilon > 0$ there exists such $N$ that for all $n > N$ 
    the following inequalities hold.
    \begin{enumerate}
        \item For all $i < r-s$ we have $$\Tilde{C}_{n,r,s}(i) < (1+\varepsilon)\Tilde{C}_{n,r,s}(i - 1).$$ 

        \item For all $i < (r-s)/2$ we have 
        $$\Tilde{C}_{n,r,s}(i) < 0.99\Tilde{C}_{n,r,s}(i - 1).$$ 
    \end{enumerate}

    Let us first notice that
    \begin{multline*}
        C_{n, r, s}(i, j) = \frac{(r-s-j+1)^2}{(r - i + 1)(n - r - i + 1)}C_{n, r, s}(i - 1, j - 1) +\\+ \frac{(s - i + j + 1)(n - 2r + s + j - i + 1)}{(r - i + 1)(n - r - i + 1)}C_{n, r, s}(i - 1, j) \leq \\ \leq \frac{(r-s-j+1)^2 + (s -i + j + 1)(n - 2r + s + j - i + 1)}{(r - i + 1)(n - r - i + 1)}\Tilde{C}_{n,r,s}(i - 1).
    \end{multline*}
    We can now see that the numerator is a quadratic function of $j$ and the coefficient before $j^2$ is positive. That is, the maximum of this function over $j = 0, \dotsc, i$ is attained either at $j = i$ or $j = 0$.

    Consider $j = i$ first. To prove the first inequality it suffices to prove that 
    $$
    (r-s-i+1)^2 + (s + 1)(n - 2r + s + 1) \leq (r - i + 1)(n - r - i + 1)
    $$
    for $i < r - s$.
    It is straightforward to check that 
    \begin{multline*}
    (r-s-i+1)^2 + (s + 1)(n - 2r + s + 1) - (r - i + 1)(n - r - i + 1) = \\ = (i - r + s) (n - 2 r + 2 s) + 1 < -n+2r-2s+1 < 0
    \end{multline*}
    and thus the first inequality is correct. The second inequality is proved in the same manner, as

    \begin{multline*}
    (r-s-i+1)^2 + (s + 1)(n - 2r + s + 1) - (r - i + 1)(n - r - i + 1) = \\ = (i - r + s) (n - 2 r + 2 s) + 1 < -\frac{1}{2}(r-s)(n-2r+2s)+1 = \\ = -\frac{1+o(1)}{2}(1-\alpha)(1-2\alpha+2\alpha^2)rn+1 \leq -\frac{1+o(1)}{8}rn < \\ < -0.01rn < -0.01(r-i+1)(n-r-i+1).
    \end{multline*}

    Let us now consider $j = 0$. Here it suffices to prove that 
    $$
    (r-s+1)^2 + (s - i + 1)(n - 2r + s - i + 1) - (r - i + 1)(n - r - i + 1) = o(n^2)
    $$
    for $i < r - s$.

    We can see that
    \begin{multline*}
        (r-s+1)^2 + (s - i + 1)(n - 2r + s - i + 1) - (r - i + 1)(n - r - i + 1) = \\ =-nr + ns + 2r^2 - 4rs + 2ir + 2s^2 - 2is + 1 = n(s-r)+2(r-s)^2+2i(r-s) + 1 = \\ = (r-s)(2(r-s)+2i-n) + 1 \leq (r-s)(4(r-s)-n) + 1.
    \end{multline*}
    It remains to notice that for $\alpha < 1/2$ we have $4(r-s)-n = (1+o(1))n(4\alpha-4\alpha^2-1) = -(1+o(1))(1-2\alpha)^2 < 0$ and for $\alpha = 1/2$ we have
    $
        4(r-s)-n = o(n).
    $

    As for the second inequality, we have 
    \begin{multline*}
        (r-s+1)^2 + (s - i + 1)(n - 2r + s - i + 1) - (r - i + 1)(n - r - i + 1) = \\ = (r-s)(2(r-s)+2i-n) + 1 \leq (r-s)(3(r-s)-n)+1 = -(r-s)(n-3r+3s)+1 = \\ = -(1+o(1))(1-3\alpha+3\alpha^2)(r-s)n+1 \leq -\frac{1+o(1)}{8}(r-s+1)n < -0.01(r-s+1)(n-r-i+1).
    \end{multline*}



    We can now proceed with the proof. Consider the following cases.
    \begin{enumerate}
        \item Suppose $i < r-s$. We can thus write
        \begin{multline*}
        |E_{r-s}(i)| = \\ = \left|\sum_{j = 0}^{i-1}(-1)^j\frac{i^2 - 2ij - i n + 2ir - 2is - 2i + jn + 2j + n s + n - r^2 - 2r + 2s + 1}{(r - i + 1)(n - r - i + 1)}\right.\times\\\times C_{n,r,s}(i-1,j)\left| \leq \frac{h(n,r,s,i) + |ns-r^2|i}{(r - i + 1)(n - r - i + 1)}\Tilde{C}_{n,r,s}(i-1)\right. = \\ = \frac{h(n,r,s,i) + |f_s(n) - 2\alpha f_r(n)|ni + f_r^2(n)i}{(r - i + 1)(n - r - i + 1)}\Tilde{C}_{n,r,s}(i-1),
    \end{multline*}
    where $h(n,r,s,i) = O(ni^2)$. We can now see that
    \begin{multline*}
        |E_{r-s}(i)| \leq \frac{h(n,r,s,i) + |f_s(n) - 2\alpha f_r(n)|ni + f_r^2(n)i}{(r - i + 1)(n - r - i + 1)}\times\\\times0.99^{\min(i-1, (r-s)/2)}1.00001^{\max(0,i-(r-s)/2)}\Tilde{C}_{n,r,s}(0) \leq  \\ \leq \frac{h(n,r,s,i) + |f_s(n) - 2\alpha f_r(n)|ni + f_r^2(n)i}{(r - i + 1)(n - r - i + 1)}0.999^{\min(i-1, (r-s)/2)}\Tilde{C}_{n,r,s}(0) \leq \\\leq \frac{h(n,r,s,i) + |f_s(n) - 2\alpha f_r(n)|ni + f_r^2(n)i}{(r - i + 1)(n - r - i + 1)}\sqrt{0.999}^{i-1}\Tilde{C}_{n,r,s}(0)
    \end{multline*}
    and thus 
    $$
        \max\limits_{i < r-s} E_{r-s}(i) = O\left(\frac{\max\{1, |f_s(n) - 2\alpha f_r(n)|, f_r^2(n)/n\}}{n}d(G(n,r,s))\right).
    $$

    \item $i \geq r-s$. Let us apply the following idea from \cite{Brouwer}. Consider the adjacency matrix $A$ of $G(n,r,s)$. Then every entry on the diagonal of $A^2$ is equal to $d = d(G(n,r,s))$ and thus $$d\binom{n}{r} = \text{tr}(A^2) = \sum_{i = 0}^{\binom{n}{r}}\lambda_i^2 = \sum_{i = 0}^{r}\left(\binom{n}{i} - \binom{n}{i - 1}\right)E_{r-s}(i)^2.$$ From this equality the following bound on $|E_{r - s}(i)|$ is easily obtained $$|E_{r - s}(i)| \leq \sqrt{\frac{\binom{n}{r}\binom{r}{s}\binom{n - r}{r - s}}{\binom{n}{i} - \binom{n}{i - 1}}} = \binom{r}{s}\binom{n - r}{r - s}\sqrt{\frac{\binom{n}{r}}{(\binom{n}{i} - \binom{n}{i - 1})\binom{r}{s}\binom{n - r}{r - s}}}.$$ 
    So if suffices to prove that $\frac{\binom{n}{r}}{\left(\binom{n}{i} - \binom{n}{i - 1}\right)\binom{r}{s}\binom{n - r}{r - s}} = O(1/n)$ for $i \geq r-s$. Note that $\binom{n}{i} - \binom{n}{i - 1} = \binom{n}{i}\frac{n-2i+1}{n-i+1} \geq \frac{1}{n}\binom{n}{r-s}$. After applying Stirling's formula, we obtain 
    \begin{multline*}
        \frac{\binom{n}{r}}{\left(\binom{n}{i} - \binom{n}{i - 1}\right)\binom{r}{s}\binom{n - r}{r - s}} \leq \frac{\binom{n}{r}}{\frac{1}{n}\binom{n}{r-s}\binom{r}{s}\binom{n - r}{r - s}} = \frac{1}{n}\frac{(r-s)!^3(n-r+s)!(n-2r+s)!s!}{(n-r)!^2r!^2} \leq \\ \leq Poly(n)\exp\{(r-s)\ln(\alpha-\alpha^2)+(n-r+s)\ln(1-\alpha+\alpha^2)\} \leq e^{-cn}.
    \end{multline*}
    \end{enumerate}   
\end{proof}

\section{Application to the giant component theorem}

In this section we highlight one possible application of our main result. To be precise, we discuss the giant component theorem. To formulate the result, let us define a random subgraph of $G(n,r,s)$ graph.

\begin{definition}
Let $G$ be a graph. Then $G(p)$ is a random graph obtained from $G$ by removing every edge with probability $1 - p$ independently from all other edges.
\end{definition}

\begin{remark}
For $G(n, r, s)$ graphs writing $G(n,r,s)(p)$ is not very convenient and thus we use the notation $G_p(n,r,s)$ instead.
\end{remark}

Let us now briefly remind the reader what is the giant component theorem. In \cite{ErdGC} the following result about the Erd\H os--R\'enyi random graph evolution was proven. It turned out that for $p < \frac{1-\varepsilon}{n}$ the graph $K_n(p)$ with high probability (w.h.p.) has only small components, while for $p > \frac{1+\varepsilon}{n}$ there w.h.p. exists a component of size $\Omega(n)$. Here we use the term $``$with high probability$"$ for such a sequence of events $\{Q_n\}_{n = 1}^{\infty}$ that $\mathbf{P}(Q_n) \to 1, n \to \infty$. It is interesting to transfer this result to graphs different from $K_n$, and indeed it can be done for various types of graphs. For example, in \cite{Yar1, Yar2, Yar3} the corresponding result was proved for $G(4s,2s,s)$ graphs and in \cite{KosEig} it was proved for $G(n,r,s)$ graphs with $n \gg r(n) \gg s(n)$. In fact, the result from \cite{KosEig} was a simple corollary of Theorem \ref{Th2} and the following general theorem from  \cite{GComp}.

\begin{theorem}
\label{GCGeneral}
Let $G_n$ be a sequence of $(n, d, \lambda)$-graphs with $\lambda = o(d)$, $n \to \infty$.
\begin{enumerate}
    \item For any $0 < \alpha < 1$ w.h.p. all components of $G_n(\frac{\alpha}{d})$ have size $O(\log(n))$.
    \item For any $\alpha > 1$ w.h.p. there exists a component of size $(1 + o(1))(1 - \frac{\overline{\alpha}}{\alpha})n$ in $G_n(\frac{\alpha}{d})$. Here, $\overline{\alpha}$ is a solution of $xe^{-x} = \alpha e^{-\alpha}$, different from $\alpha$. At the same time, all other components of $G_n(\frac{\alpha}{d})$ have size $O(\log(n))$.
\end{enumerate}
\end{theorem}

Using Theorem \ref{Main} we can obtain similar result for $G(n, \alpha n + f_r(n), \alpha^2 n + f_s(n))$ graphs.

\begin{theorem}
\label{GC1aa2}
Let $G_n = G(n, \alpha n + f_r(n), \alpha^2n+f_s(n))$, $G_n(p_n) = G(n, \alpha n + f_r(n), \alpha^2n+f_s(n))$, $d_n = d(G_n)$, $N_n = \binom{n}{\alpha n + f_r(n)}$. Assume that $f_r(n) = o(n), f_s(n) = o(n)$. Then
\begin{enumerate}
    \item For any $0 < \alpha < 1$ w.h.p. all components of $G_n(\frac{\alpha}{d_n})$ have size $O(\log(N_n))$.
    \item For any $\alpha > 1$ w.h.p. there exists a component of size $(1 + o(1))(1 - \frac{\overline{\alpha}}{\alpha})N_n$ in $G_s(\frac{\alpha}{d_n})$. Here, $\overline{\alpha}$ is a solution of $xe^{-x} = \alpha e^{-\alpha}$, different from $\alpha$. At the same time, all other components of $G_n(\frac{\alpha}{d_n})$ have size $O(\log(N_n))$.
\end{enumerate}
\end{theorem}

\begin{proof}
    Indeed, as both $|f_s(n) - 2\alpha f_r(n)|$ and $f_r^2(n)/n$ are $o(n)$, the assumption of \ref{GCGeneral} holds.
\end{proof}

\end{document}